\documentclass[amscd,amssymb,11pt]{amsart}

\usepackage{amsmath,amssymb}
\usepackage[all]{xy}

\theoremstyle{definition}

\theoremstyle{remark}

\numberwithin{equation}{section}

\newcommand{\Z}{{\mathbb  Z}}

\newcommand{\R}{{\mathbb  R}}

\newcommand{\Sinfty}{\Sigma^{\infty}}

\newcommand{\ra}{\rightarrow}

\begin{document}

\title[Correction to Nonrealization results via Goodwillie towers]{Correction to ``Topological nonrealization results via the Goodwillie tower approach to iterated loopspace homology''}

\author[Kuhn]{Nicholas J.~Kuhn}
\address{Department of Mathematics \\ University of Virginia \\ Charlottesville, VA 22904}
\email{njk4x@virginia.edu}
\thanks{This research was partially supported by grants from the National Science Foundation}

\date{November 9, 2009.}


\begin{abstract}  Manfred Stelzer has pointed out that part of Corollary 4.5 of \cite{k2} was not sufficiently proved, and, indeed, is likely incorrect as stated.  This necessitates a little more argument to finish the proof of the main theorem of \cite{k2}.  The statement of this theorem, and all the examples, remain unchanged.

\end{abstract}

\maketitle

In \cite{k2}, the author showed that certain unstable modules over the mod 2 Steenrod algebra couldn't be realized as the reduced mod 2 cohomology of a space.  The modules have the form $\Sigma^n M$, where $M$ is an unstable module of a special sort.  The method of proof was to use a 2nd quadrant spectral sequence converging to $H^*(\Omega^n X; \Z/2)$ to show that, were a space $X$ to exist whose cohomology realized $\Sigma^n M$, $H^*(\Omega^n X;\Z/2)$ could not admit a cup product compatible with Steenrod operations.

The spectral sequence for $n >1$ is a newish one, arising from the Goodwillie tower of the functor $X \mapsto \Sinfty \Omega^n X$, and section 4 of \cite{k2} is devoted to collecting and proving some basic facts about this spectral sequence.  I thank Manfred Stelzer for pointing out that part of Corollary 4.5 is likely over optimistic, and certainly was not sufficiently proved.

We assume notation as in \cite{k2}.

In Corollary 4.5, it was asserted that, if $\tilde H^*(X; \Z/2) \simeq \Sigma^n M$ with $M$ unstable, and also has no nontrivial cup products, then in the spectral sequence, one will have $E_3^{-1,*} = E_2^{-1,*} = E_1^{-1,*}$, and $E_2^{-2,*} = E_1^{-2,*}$.  My mistake was in not adequately considering possible differentials on elements in $E_1^{-3,*}$ of the form $\sigma^3 L_{n-1}(x\otimes y \otimes z)$.  Under the hypotheses on the cup product, the $d_1$ differential on such terms {\em will} be 0, by the same argument given explaining why $d_1$ is zero on terms of the form $L(x \otimes y)$: by comparison to the classical Eilenberg--Moore spectral sequence.  But there is no apparent reason why $d_2$ need also be zero on such terms.  We can only conclude that $E_2^{-1,*} = E_1^{-1,*}$, and $E_2^{-2,*} = E_1^{-2,*}$.

Corollary 4.5 is used at one critical point in the proof of the main theorem given in section 5.  Lemma 5.3 asserts that a certain element in $E_1^{-1,2d+2^{k+2}+1}$ is not a boundary.  The argument given is that for dimension reasons, no $d_r$ for $r >2$ could have nonzero image in this bigrading.  Implicit is that Corollary 4.5 takes care of $d_1$ and $d_2$. In light of the comments above, one needs a new argument for $d_2$.

It turns out that, except for one special case, a dimension argument still works: $E_3^{-3, 2d+ 2^{k+2} + 2}$ contains no elements of the form $\sigma^3 L_{n-1}(x \otimes y \otimes z)$.  There are two extreme cases to consider: if $x$, $y$, and $z$ are all chosen from the top of $N_0$, and if $x$ and $y$ are chosen from the bottom of $N_0$ and $z$ is chosen from the bottom of $M_1$.

In the first case, $|x|=|y|=|z| = m + 2^k$, and so $\sigma^3 L_{n-1}(x \otimes y \otimes z)$ has bidegree $(-3, 3m+ 3 \cdot 2^{k} + 2n+1)$.  In the second case, $|x|=|y| = d+2^k$ and $|z| = l+2^{k+1}$, and so $\sigma^3 L_{n-1}(x \otimes y \otimes z)$ has bidegree $(-3, 2d+ l + 2^{k+2} + 2n+1)$.

We are assuming inequality (5--3), which says that $2^k > 4m-2l+2n-2$.  One also has that $0 \leq l \leq d \leq m$ and $n \geq 1$.  One can then check that, indeed,
$$ 3m+ 3 \cdot 2^{k} + 2n+1 < 2d+ 2^{k+2} + 2 < 2d+ l + 2^{k+2} + 2n+1,$$
{\em unless} we are in the special case $k=0, n=1, l=d=m=0$.

In this final special case, $n=1$, so we are trying to use the classical Eilenberg--Moore spectral sequence to show that, if $M$ is a $\Z/2$ vector space concentrated in degree 0, there cannot exist a space $X$ with $\tilde H^*(X;\Z/2) \simeq \Sigma M \otimes \Phi(0,2)$, if all cup products are zero.  Such a space will necessarily fit into a cofibration sequence of the form
$$ \bigvee S^4 \ra \bigvee \Sigma \R P^2 \ra X.$$
We leave it to the reader to check that, by appropriately including $S^4$ into the first wedge, and projecting out onto a $\Sigma \R P^2$ in the second wedge, one sees that $X$ will have a `subquotient' $Y$ with
$\tilde H^*(Y;\Z/2) \simeq \Sigma \Phi(0,2)$, and still with all cup products 0.

Similar to, but simpler than, arguments in section 6 of \cite{k2} (which dealt with $\Sigma^2 \Phi(1,3)$), our arguments show that such a $Y$ can't exist.  Repressing some suspensions from the notation, Figure 1 shows all of $E_1^{*,*}$ in total degree less than or equal to 4, in the Eilenberg--Moore spectral sequence converging to $H^*(\Omega Y;\Z/2)$.

\begin{figure}
\begin{equation*}
\begin{array}{ccccc|c}
a \otimes a \otimes a \otimes a& &&&& 8\\
& & &&& 7\\
 &a \otimes a \otimes a &b \otimes b&&& 6\\
&   &a\otimes b, b\otimes a&c&& 5\\
 &   &a\otimes a&&& 4\\
 &    &&b&& 3\\
&     &&a&& 2\\
& &&&& 1\\
& &&&  1   & 0  \\ \hline
-4& -3&-2&-1& 0&s\backslash t \\
\end{array}
\end{equation*}
\caption{$E_1^{s,t}$ when $\tilde H^*(Y;\Z/2) \simeq \Sigma\Phi(0,2)$} \label{figure 1}
\end{figure}

As cup products are assumed zero, $E_2^{*,*}= E_1^{*,*}$.  Furthermore, $d_2(a\otimes a \otimes a) = 0$ (and thus {\em not} $c$), because $a \otimes a \otimes a = (a \otimes a) * a$ and $d_2$ is a derivation with respect to the shuffle product $*$. Thus through degree 4, $F^{-2}H^*(\Omega Y;\Z/2)$ would have a basis given by elements 1, $\alpha$, $\beta$, $\delta$, $\epsilon_1$, $\epsilon_2$, $\gamma$, and $\omega$, in respective degrees 0, 1, 2, 2, 3, 3, 4, and 4, and represented by $1$, $a$, $b$, $a \otimes a$, $a \otimes b$, $b \otimes a$, $c$, and $b \otimes b$.  The structure of $\Phi(0,2)$ ($Sq^1 a = b$, $Sq^2 b = c$) shows that $\gamma = \beta^2 = \alpha^4$.  Furthermore, $Sq^1 \delta = \epsilon_1 + \epsilon_2 = \alpha \cup \beta$, as all three are represented by $a \otimes b + b \otimes a$.  One then gets a contradiction, as
$$ 0 = Sq^1 Sq^1 \delta = Sq^1(\alpha \cup \beta) = \beta^2 = \gamma \neq 0.$$

We end by observing that $\tilde H^*(SU(3)/SO(3);\Z/2) \simeq \Sigma \Phi(0,2)$. Here, of course, cup products are not zero, due to Poincar\'e duality.

\end{document}